\documentclass[12pt,fleqn,leqno]{article}

\author{
D.D.\ Poro\c sniuc\\  }
\date{}
\title{A locally symmetric K\"ahler Einstein structure on a tube in the nonzero cotangent
bundle of a space form}

\begin{document}

\maketitle
\begin{abstract}
We obtain a locally symmetric K\"ahler Einstein structure on a
tube in the nonzero cotangent bundle of a Riemannian manifold of
positive constant sectional curvature. The obtained K\"ahler
Einstein structure cannot have constant holomorphic sectional
curvature.


\end{abstract}

\vskip5mm {\large \bf Introduction} \vskip5mm The differential
geometry of the cotangent bundle $T^*M$ of a Riemannian manifold
$(M,g)$ is almost similar to that of the tangent bundle $TM$.
However, there are some differences because the lifts (vertical,
complete, horizontal etc.) to $T^*M$ cannot be defined just like
in the case of $TM$.

In \cite{OprPor} V. Oproiu and the present author have obtained a
natural K\"ahler Einstein structure $(G,J)$ of diagonal type
induced on $T^*M$ from the Riemannian metric $g$. The obtained
K\"ahler structure on $T^*M$ depends on one essential parameter
$u$, which is a smooth function depending on the energy density
$t$ on $T^*M$. If the K\"ahler structure is Einstein they get a
second order differential equation fulfilled by the parameter $u$.
In the case of the general solution, they have obtained that
$(T^*M,G,J)$ has constant holomorphic sectional curvature.

 In this paper we study the singular case where the parameter
$u=At$, $A \in {\bf R}$. The considered natural Riemannian metric
$G$ of diagonal type on the nonzero cotangent bundle $T^*_0M$ is
defined by using one parameter $v$ which is a smooth function
depending on the energy density $t$. The vertical distribution
$VT^*_0M$ and the horizontal distribution $HT^*_0M$ are orthogonal
to each other but the dot products induced on them from $G$ are
not isomorphic (isometric).

 Next, the natural almost complex structures $J$ on
$T^*_0M$ that interchange the vertical and horizontal
distributions depends of one essential parameter $v$.

 After that, we obtain that $G$ is Hermitian with respect to $J$ and it
 follows that the fundamental $2$-form $\phi$ associated to the
 almost Hermitian structure $(G,J)$ is the fundamental form defining the
 usual symplectic structure on $T^*_0M$, hence it is closed.

  From the integrability condition for $J$ it
follows that the base manifold $M$ must have constant sectional
curvature $c$ and the parameter $v$ must be a rational function
depending on energy density $t$.

 If the constant sectional curvature c is positive then
we obtain a locally symmetric K\"ahler Einstein structure defined
on a tube in $T^*_0M$. The K\"ahler Einstein manifold obtained
cannot have constant holomorphic sectional curvature.

 The manifolds, tensor fields and geometric objects we consider
in this paper, are assumed to be differentiable of class
$C^{\infty}$ (i.e. smooth). We use the computations in local
coordinates but many results from this paper may be expressed in
an invariant form. The well known summation convention is used
throughout this paper, the range for the indices $h,i,j,k,l,r,s$
being always${\{}1,...,n{\}}$ (see \cite{GheOpr}, \cite{OprPap1},
\cite{OprPap2}, \cite{YanoIsh}). We shall denote by
${\Gamma}(T^*_0M)$ the module of smooth vector fields on $T^*_0M$.

\vskip5mm {\large \bf 1. Some geometric properties of $T^*M$}
\vskip5mm

Let $(M,g)$ be a smooth $n$-dimensional Riemannian manifold and
denote its cotangent bundle by $\pi :T^*M\longrightarrow M$.
Recall that there is a structure of a $2n$-dimensional smooth
manifold on $T^*M$, induced from the structure of smooth
$n$-dimensional manifold  of $M$. From every local chart
$(U,\varphi )=(U,x^1,\dots ,x^n)$  on $M$, it is induced a local
chart $(\pi^{-1}(U),\Phi )=(\pi^{-1}(U),q^1,\dots , q^n,$
$p_1,\dots ,p_n)$, on $T^*M$, as follows. For a cotangent vector
$p\in \pi^{-1}(U)\subset T^*M$, the first $n$ local coordinates
$q^1,\dots ,q^n$ are  the local coordinates $x^1,\dots ,x^n$ of
its base point $x=\pi (p)$ in the local chart $(U,\varphi )$ (in
fact we have $q^i=\pi ^* x^i=x^i\circ \pi, \ i=1,\dots n)$. The
last $n$ local coordinates $p_1,\dots ,p_n$ of $p\in \pi^{-1}(U)$
are the vector space coordinates of $p$ with respect to the
natural basis $(dx^1_{\pi(p)},\dots , dx^n_{\pi(p)})$, defined by
the local chart $(U,\varphi )$,\ i.e. $p=p_idx^i_{\pi(p)}$.

An $M$-tensor field of type $(r,s)$ on $T^*M$ is defined by sets
of $n^{r+s}$ components (functions depending on $q^i$ and $p_i$),
with $r$ upper indices and $s$ lower indices, assigned to induced
local charts $(\pi^{-1}(U),\Phi )$ on $T^*M$, such that the local
coordinate change rule is that of the local coordinate components
of a tensor field of type $(r,s)$ on the base manifold $M$ (see
\cite{Mok} for further details in the case of the tangent bundle).
An usual tensor field of type $(r,s)$ on $M$ may be thought of as
an $M$-tensor field of type $(r,s)$ on $T^*M$. If the considered
tensor field on $M$ is covariant only, the corresponding
$M$-tensor field on $T^*M$ may be identified with the induced
(pullback by $\pi $) tensor field on $T^*M$.

Some useful $M$-tensor fields on $T^*M$ may be obtained as
follows. Let $v,w:[0,\infty ) \longrightarrow {\bf R}$ be a smooth
functions and let $\|p\|^2=g^{-1}_{\pi(p)}(p,p)$ be the square of
the norm of the cotangent vector $p\in \pi^{-1}(U)$ ($g^{-1}$ is
the tensor field of type (2,0) having the components $(g^{kl}(x))$
which are the entries of the inverse of the matrix $(g_{ij}(x))$
defined by the components of $g$ in the local chart $(U,\varphi
)$). The components $g_{ij}(\pi(p))$, $p_i$, $v(\|p\|^2)p_ip_j $
define $M$-tensor fields of types $(0,2)$, $(0,1)$, $(0,2)$ on
$T^*M$, respectively. Similarly, the components $g^{kl}(\pi(p))$,
$g^{0i}=p_hg^{hi}$, $w(\|p\|^2)g^{0k}g^{0l}$ define $M$-tensor
fields of type $(2,0)$, $(1,0)$, $(2,0)$ on $T^*M$, respectively.
Of course, all the components considered above are in the induced
local chart $(\pi^{-1}(U),\Phi)$.

 The Levi Civita connection $\dot \nabla $ of $g$ defines a direct
sum decomposition

\begin{equation}
TT^*M=VT^*M\oplus HT^*M.
\end{equation}
of the tangent bundle to $T^*M$ into vertical distributions
$VT^*M= {\rm Ker}\ \pi _*$ and the horizontal distribution
$HT^*M$.

 If $(\pi^{-1}(U),\Phi)=(\pi^{-1}(U),q^1,\dots ,q^n,p_1,\dots ,p_n)$
is a local chart on $T^*M$, induced from the local chart
$(U,\varphi )= (U,x^1,\dots ,x^n)$, the local vector fields
$\frac{\partial}{\partial p_1}, \dots , \frac{\partial}{\partial
p_n}$ on $\pi^{-1}(U)$ define a local frame for $VT^*M$ over $\pi
^{-1}(U)$ and the local vector fields $\frac{\delta}{\delta
q^1},\dots ,\frac{\delta}{\delta q^n}$ define a local frame for
$HT^*M$ over $\pi^{-1}(U)$, where
$$
\frac{\delta}{\delta q^i}=\frac{\partial}{\partial
q^i}+\Gamma^0_{ih} \frac{\partial}{\partial p_h},\ \ \ \Gamma
^0_{ih}=p_k\Gamma ^k_{ih}
 $$
and $\Gamma ^k_{ih}(\pi(p))$ are the Christoffel symbols of $g$.

The set of vector fields $(\frac{\partial}{\partial p_1},\dots
,\frac{\partial}{\partial p_n}, \frac{\delta}{\delta q^1},\dots
,\frac{\delta}{\delta q^n})$ defines a local frame on $T^*M$,
adapted to the direct sum decomposition (1).

We consider
\begin{equation}
t=\frac{1}{2}\|p\|^2=\frac{1}{2}g^{-1}_{\pi(p)}(p,p)=\frac{1}{2}g^{ik}(x)p_ip_k,
\ \ \ p\in \pi^{-1}(U)
\end{equation}
the energy density defined by $g$ in the cotangent vector $p$. We
have $t\in [0,\infty)$ for all $p\in T^*M$.

From now on we shall work in a fixed local chart $(U,\varphi)$ on
$M$ and in the induced local chart $(\pi^{-1}(U),\Phi)$ on $T^*M$.

\vskip5mm {\large \bf 2. An almost K\"ahler structure on the
$T^*_0M$} \vskip5mm

The nonzero cotangent bundle $T^*_0M$ of Riemannian manifold
$(M,g)$ is defined by the formula: $T^*M$ minus zero section.
 Consider a real valued smooth
function $v$ defined on $(0,\infty)\subset {\bf R}$ and a real
constant A. We define the following $M$-tensor field of type
$(0,2)$  on $T^*_0M$ having the components
\begin{equation}
G_{ij}(p)= At g_{ij}(\pi(p))+v(t)p_ip_j.
\end{equation}

It follows easily that the matrix $(G_{ij})$ is positive definite
if and only if $A>0,\ A+2v>0$. The inverse of this matrix has the
entries

\begin{equation}
H^{kl}(p)= \frac{1}{At} g^{kl}(\pi(p))+w(t)g^{0k}g^{0l},
\end{equation}
where
\begin{equation}
w= \frac{-v}{At^2(A+2v)}.
\end{equation}

The components $H^{kl}$ define an $M$-tensor field of type $(2,0)$
on $T^*_0M$ .\\

{\bf Remark.} If the matrix $(G_{ij})$ is positive definite, then
its inverse $(H^{kl})$ is positive definite too.

Using the $M$-tensor fields defined by $G_{ij},\ H^{kl}$, the
following Riemannian metric may be considered on $T^*_0M$
\begin{equation}
G=G_{ij}dq^idq^j+H^{ij}Dp_iDp_j,
\end{equation}
where $Dp_i=dp_i-\Gamma^0_{ij}dq^j$ is the absolute (covariant)
differential of $p_i$ with respect to the Levi Civita connection
$\dot\nabla$ of $g$ . Equivalently, we have
$$
G(\frac{\delta}{\delta q^i},\frac{\delta}{\delta
q^j})=G_{ij},~~~G(\frac{\partial}{\partial p_i}
,\frac{\partial}{\partial p_j})=H^{ij},~~
G(\frac{\partial}{\partial p_i},\frac{\delta}{\delta q^j})=~
G(\frac{\delta}{\delta q^j},\frac{\partial}{\partial p_i})=0.
$$

Remark that $HT^*_0M,~VT^*_0M$ are orthogonal to each other with
respect to $G$, but the Riemannian metrics induced from $G$ on
$HT^*_0M,~VT^*_0M$ are not the same, so the considered metric $G$
on $T^*M$ is not a metric of Sasaki type.
 Remark also that the system of 1-forms
$(dq^1,...,dq^n,Dp_1,...,Dp_n)$ defines a local frame on
$T^{*}T^*M$, dual to the local frame $(\frac{\delta}{\delta
q^1},\dots ,\frac{\delta}{\delta q^n},\frac{\partial}{\partial
p_1},\dots ,\frac{\partial}{\partial p_n})$ adapted to the direct
sum decomposition (1).

Next, an almost complex structure $J$ is defined on $T^*_0M$ by
the same $M$-tensor fields $G_{ij},\ H^{kl}$, expressed in the
adapted local frame by
\begin{equation}
J\frac{\delta}{\delta
q^i}=G_{ik}\frac{\partial}{\partial p_k},\ \ \
J\frac{\partial}{\partial p_i}=-H^{ik}\frac{\delta}{\delta q^k}
\end{equation}

 From the property
of the $M$-tensor field $H^{kl}$ to be defined by the inverse of
the matrix defined by the components of the $M$-tensor field
$G_{ij}$, it follows easily that $J$ is an almost complex
structure on $T^*_0M$. \\

{\bf Theorem 1}. \it $(T^*_0M,G,J)$ is an almost K\"ahler
manifold.\\

 Proof. \rm Since the matrix $(H^{kl})$ is the inverse of the
 matrix $(G_{ij})$, it follows easily that
$$
G(J\frac{\delta}{\delta q^i},J\frac{\delta}{\delta
q^j})=G(\frac{\delta}{\delta q^i},\frac{\delta}{\delta q^j}),\ \ \
G(J\frac{\partial}{\partial p_i},J\frac{\partial}{\partial
p_j})=G(\frac{\partial}{\partial p_i},\frac{\partial}{\partial
p_j}),$$ $$ G(J\frac{\partial}{\partial p_i},J\frac{\delta}{\delta
q^j})=G(\frac{\partial}{\partial p_i},\frac{\delta}{\delta
q^j})=0.
$$

Hence
$$G(JX,JY)=G(X,Y),\ \ \forall\ X,Y{\in}{\Gamma}(T^*_0M).$$
Thus  $(T^*_0M,G,J)$ is an almost Hermitian manifold.

The fundamental $2$-form associated with this almost Hermitian
structure is $\phi$, defined by
$$\phi(X,Y) = G(X,JY),\ \ \ \forall\ X,Y{\in}{\Gamma}(T^*_0M).$$
By a straightforward computation we get
$$
\phi(\frac{\delta}{\delta q^i},\frac{\delta}{\delta q^j})=0,\ \ \
\phi(\frac{\partial}{\partial p_i},\frac{\partial}{\partial
p_j})=0,\ \ \ \phi(\frac{\partial}{\partial
p_i},\frac{\delta}{\delta q^j})= \delta^i_j.
$$
Hence
\begin{equation}
\phi =Dp_i\wedge dq^i= dp_i\wedge dq^i,
\end{equation}
due to the symmetry of $\Gamma^0_{ij}=p_h\Gamma^h_{ij}$. It
follows that $\phi$ does coincide with the fundamental $2$-form
defining the usual symplectic structure on $T^*_0M$. Of course, we
have $d\phi =0$, i.e. $\phi$ is closed. Therefore $(T^*_0M,G,J)$
is an almost K\"ahler manifold.

\vskip5mm {\large \bf 3. A K\"ahler structure on a tube
$T^*_{0A}M$} \vskip5mm

We shall study the integrability of the almost complex structure
defined by $J$ on $T^*_0M$. To do this we need the following well
known formulas for the brackets of the vector fields
$\frac{\partial}{\partial p_i},\frac{\delta}{\delta q^i},~
i=1,...,n$
\begin{equation}
[\frac{\partial}{\partial p_i},\frac{\partial}{\partial
p_j}]=0;~~~[\frac{\partial}{\partial p_i},\frac{\delta}{\delta
q^j}]=\Gamma^i_{jk}\frac{\partial}{\partial p_k};~~~
[\frac{\delta}{\delta q^i},\frac{\delta}{\delta q^j}]
=R^0_{kij}\frac{\partial}{\partial p_k},
\end{equation}
where $R^h_{kij}(\pi(p))$ are the local coordinate components of
the curvature tensor field of $\dot \nabla$ on $M$ and
$R^0_{kij}(p)=p_hR^h_{kij}$. Of course, the components
 $R^0_{kij}$, $R^h_{kij}$ define M-tensor fields of types
 (0,3), (1,3) on $T^*M$, respectively.\\

{\bf Theorem 2. } {\it The Nijenhuis tensor field of the almost
complex structure $J$ on $T^*_0M$ is given by}
\begin{equation}
\left\{
\begin{array}{l}
N(\frac{\delta}{\delta q^i},\frac{\delta}{\delta
q^j})=\{At(v+A)(\delta^h_ig_{jk}-
\delta^h_jg_{ik})-R^h_{kij}\}p_h\frac{\partial}{\partial p_k},
\\ \mbox{ } \\
N(\frac{\delta}{\delta q^i},\frac{\partial}{\partial
p_j})=H^{kl}H^{jr}\{At(v+A)(\delta^h_ig_{rl}-
\delta^h_rg_{il})-R^h_{lir}\}p_h\frac{\delta}{\delta q^k},
\\ \mbox{ } \\
N(\frac{\partial}{\partial p_i},\frac{\partial}{\partial
p_j})=H^{ir}H^{jl}\{At(v+A)(\delta^h_lg_{rk}-
\delta^h_rg_{lk})-R^h_{klr}\}p_h\frac{\partial}{\partial p_k}.
\end{array}
\right.
\end{equation}

\vskip2mm

{\it Proof. } Recall that the Nijenhuis tensor field $N$ defined
by $J$ is given by
$$
N(X,Y)=[JX,JY]-J[JX,Y]-J[X,JY]-[X,Y],\ \ \forall\ \ X,Y \in \Gamma
(T^*_0M).
$$
Then, we have $\frac{\delta}{\delta q^k}t =0,\
\frac{\partial}{\partial p_k}t = g^{0k}$ and $\dot
\nabla_iG_{jk}=0,\dot \nabla_iH^{jk}= 0$, where
$$
\dot \nabla_iG_{jk}= \frac{\delta}{\delta
q^i}G_{jk}-\Gamma^l_{ij}G_{lk}-\Gamma^l_{ik}G_{lj}
$$
$$
\dot \nabla_iH^{jk}= \frac{\delta}{\delta
q^i}H^{jk}+\Gamma^j_{il}H^{lk}+\Gamma^k_{il}H^{lj}
$$

The above expressions for the components of $N$ can be obtained by
a quite long, straightforward  computation.\\

{\bf Theorem 3.} {\it The almost complex structure $J$ on $T^*_0M$
is integrable if and only if the base manifold $M$ has constant
sectional curvature $c$ and the function $v$ is given by}
\begin{equation}
v=\frac {c-A^2t}{At}.
\end{equation}

{\it Proof.} From the condition $N=0$, one obtains
$$
\{At(v+A)(\delta^h_ig_{jk}- \delta^h_jg_{ik})-R^h_{kij}\}p_h=0.
$$
Taking $p_l\neq 0$ and $p_h=0~\forall h\neq l$, it follows that
$$
R^l_{kij}=At(v+A)(\delta^l_ig_{jk}- \delta^l_jg_{ik}).
$$
 Thus the sectional curvature $c=At(v+A)$ of $(M,g)$ depends
 only on $q^i$. Using the Schur theorem (in the case where $M$
 is connected and ${\rm dim}\ M \geq 3$), we obtain that c is constant .
 Then we obtain the expression (11)
 of $v$.

Conversely, if $(M,g)$ has constant sectional curvature $c$ and
$v$ is given by (11), it
follows in a straightforward way that $N = 0$.\\

\bf Remark. \rm  The function $v$ must fulfill the condition
\begin{equation}
A+2v=\frac{2c-A^2t}{At}>0,~~A>0.
\end{equation}

If $c>0$ then $(T^*_{0A}M,G,J)$ is a K\"ahler manifold, where
$T^*_{0A}M$ is the tube in $T^*_0M$ defined by the condition
$0<\|p\|^2<\frac{4c}{A^2}.$

 The components of the K\"ahler metric $G$ on $T^*_{0A}M$ are

\begin{equation}
\left\{
\begin{array}{l}
G_{ij}=Atg_{ij}+\frac
{c-A^2t}{At}p_ip_j,
\\ \mbox{ } \\
H^{ij}=\frac{1}{At}g^{ij}-\frac{c-A^2t}{At^2(2c-A^2t)}g^{0i}g^{0j}.
\end{array}
\right.
\end{equation}
\newpage
\vskip5mm
 {\large \bf 4. A  K\"ahler  Einstein structure on
 $T^*_{0A}M$}
\vskip5mm
In this section we shall study the property of the
K\"ahler manifold $(T^*_{0A}M,G,J)$ to be Einstein.

The Levi Civita connection $\nabla$ of the Riemannian manifold
$(T^*_{0A}M,G)$ is determined by the conditions
$$
\nabla G=0,~~~~~  T =0,
$$
where $T$ is its torsion tensor field. The explicit expression of
this connection is obtained from the formula
$$
2G({\nabla}_XY,Z)=X(G(Y,Z))+Y(G(X,Z))-Z(G(X,Y))+
$$
$$
+G([X,Y],Z)-G([X,Z],Y)-G([Y,Z],X); ~~~~~~ \forall\ X,Y,Z~{\in}~
{\Gamma}(T^*_{0A}M).
$$

The final result can be stated as follows. \\

\bf Theorem 4. {\it The Levi Civita connection ${\nabla}$ of $G$
has the following expression in the local adapted frame
$(\frac{\delta}{\delta q^1},\dots ,\frac{\delta}{\delta
q^n},\frac{\partial}{\partial p_1},\dots ,\frac{\partial}{\partial
p_n}):$
\begin{equation}
\left\{
\begin{array}{l}
 \nabla_\frac{\partial}{\partial
p_i}\frac{\partial}{\partial p_j}
=Q^{ij}_h\frac{\partial}{\partial p_h},\ \ \ \ \ \
\nabla_\frac{\delta}{\delta q^i}\frac{\partial}{\partial
p_j}=-\Gamma^j_{ih}\frac{\partial}{\partial
p_h}+P^{hj}_i\frac{\delta}{\delta q^h},
\\ \mbox{ } \\
\nabla_\frac{\partial}{\partial p_i}\frac{\delta}{\delta
q^j}=P^{hi}_j\frac{\delta}{\delta q^h},\ \ \ \ \ \
\nabla_\frac{\delta}{\delta q^i}\frac{\delta}{\delta
q^j}=\Gamma^h_{ij}\frac{\delta}{\delta
q^h}+S_{hij}\frac{\partial}{\partial p_h},
\end{array}
\right.
\end{equation}
where $Q^{ij}_h, P^{hi}_j, S_{hij}$ are $M$-tensor fields on
$T^*_{0A}M$, defined by
\begin{equation}
\left\{
\begin{array}{l}
 Q^{ij}_h = \frac{1}{2}G_{hk}(\frac{\partial}{\partial
p_i}H^{jk}+ \frac{\partial}{\partial p_j}H^{ik}
-\frac{\partial}{\partial p_k}H^{ij}),
\\ \mbox{ } \\
 P^{hi}_j=\frac{1}{2}H^{hk}(\frac{\partial}{\partial
p_i}G_{jk}-H^{il}R^0_{ljk}),
\\ \mbox{ } \\
S_{hij}=-\frac{1}{2}G_{hk}\frac{\partial}{\partial
p_k}G_{ij}+\frac{1}{2}R^0_{hij}.
\end{array}
\right.
\end{equation}
\rm

After replacing of the expressions of the involved $M$-tensor
fields , we obtain

\begin{equation}
\left\{
\begin{array}{l}
Q^{ij}_h =
\frac{1}{2t}[(g^{ij}+\frac{c}{t(2c-A^2t)}g^{0i}g^{0j})p_h-(\delta^i_hg^{0j}+\delta^j_hg^{0i})],
\\ \mbox{ } \\
P^{hi}_j=-Q^{ih}_j,
\\ \mbox{ } \\
S_{hij}=\frac{-2c+A^2t}{2}(g_{ij}p_h+g_{ih}p_j)+\frac{A^2t}{2}g_{hj}p_{i}+\frac{3c-2A^2t}{2t}p_hp_ip_j
.
\end{array}
\right.
\end{equation}

The curvature tensor field $K$ of the connection $\nabla $ is
obtained from the well known formula
$$
K(X,Y)Z=\nabla_X\nabla_YZ-\nabla_Y\nabla_XZ-\nabla_{[X,Y]}Z,\ \ \
\ \forall\ X,Y,Z\in \Gamma (T^*_{0A}M).
$$

The components of curvature tensor field $K$ with respect to the
adapted local frame $(\frac{\delta}{\delta q^1},\dots
,\frac{\delta}{\delta q^n},\frac{\partial}{\partial p_1},\dots
,\frac{\partial}{\partial p_n})$ are obtained easily:
\begin{equation}
\left\{
\begin{array}{l}
K(\frac{\delta}{\delta q^i},\frac{\delta}{\delta
q^j})\frac{\delta}{\delta q^k}=QQQ^h_{ijk}\frac{\delta}{\delta
q^h},\ \ \ \ \ K(\frac{\delta}{\delta q^i},\frac{\delta}{\delta
q^j})\frac{\partial}{\partial
p_k}=QQP^k_{ijh}\frac{\partial}{\partial p_h},
\\ \mbox{ } \\
K(\frac{\partial}{\partial p_i},\frac{\partial}{\partial
p_j})\frac{\delta}{\delta q^k}=PPQ^{ijh}_k\frac{\delta}{\delta
q^h},\ \ \ \ \ K(\frac{\partial}{\partial
p_i},\frac{\partial}{\partial p_j})\frac{\partial}{\partial p_k}
=PPP^{ijk}_h\frac{\partial}{\partial p_h},
\\ \mbox{ } \\
K(\frac{\partial}{\partial p_i},\frac{\delta}{\delta
q^j})\frac{\delta}{\delta q^k}=PQQ^i_{jkh}\frac{\partial}{\partial
p_h},\ \ \ \ \ K(\frac{\partial}{\partial
p_i},\frac{\delta}{\delta q^j})\frac{\partial}{\partial
p_k}=PQP^{ikh}_j\frac{\delta}{\delta q^h},
\end{array}
\right.
\end{equation}
where
\begin{equation}
\left\{
\begin{array}{l}
QQQ^h_{ijk}=\frac{A^2t}{2}(\delta^h_ig_{jk}-\delta^h_jg_{ik})+\frac{A^2}{4}(g_{ik}p_j-g_{jk}p_i)g^{0h}-
\\ \mbox{ } \\
 ~~~~~~~~~~~~~-\frac{A^2}{4}(\delta^h_ip_j-\delta^h_jp_i)p_k,
\\ \mbox{ } \\
QQP^k_{ijh}=-QQQ^k_{ijh},
\\ \mbox{ } \\
PPQ^{ijh}_k=-\frac{1}{2t}(\delta^i_kg^{jh}-\delta^j_kg^{ih})-\frac{1}{4t^2}(g^{ih}g^{0j}-g^{jh}g^{0i})p_k+
\\ \mbox{ } \\
~~~~~~~~~~~~~+\frac{1}{4t^2}(\delta^i_kg^{0j}-\delta^j_kg^{0i})g^{0h},
\\ \mbox{ } \\
PPP^{ijk}_h=-PPQ^{ijk}_h,
\\ \mbox{ } \\
PQQ^i_{jkh}=\frac{A}{2}\delta^i_jG_{hk}+\frac{2c-A^2t}{4t}(\delta^i_kp_h+\delta^i_hp_k)p_j+
\\ \mbox{ } \\
~~~~~~~~~~~~~+\frac{A^2}{4}(g_{jh}p_k+g_{jk}p_h)g^{0i}-\frac{c}{2t^2}g^{0i}p_jp_hp_k,
\\ \mbox{ } \\
PQP^{ikh}_j=-\frac{A}{2}\delta^i_jH^{hk}-\frac{1}{4t^2}(g^{ih}g^{0k}+g^{ik}g^{0h})p_j-
\\ \mbox{ } \\
~~~~~~~~~~~~~~-\frac{A^2}{4t(2c-A^2t)}(\delta^h_jg^{0k}+\delta^k_jg^{0h})g^{0i}
+\frac{c}{2t^3(2c-A^2t)}g^{0i}g^{0h}g^{0k}p_j,
\end{array}
\right.
\end{equation}
are M-tensor fields on $T^*_{0A}M$.\\

 {\bf Remark.} From the local coordinates expression of the
curvature tensor field K we obtain that the K\"ahler manifold
$(T^*M,G,J)$ cannot have constant holomorphic sectional curvature.\\

The Ricci tensor field Ric of $\nabla$ is defined by the formula:
$$
Ric(Y,Z)=trace(X\longrightarrow K(X,Y)Z),\ \ \ \forall\  X,Y,Z\in
\Gamma (T^*_{0A}M).
$$
It follows
$$
\left\{
\begin{array}{l}
Ric(\frac{\delta}{\delta q^i},\frac{\delta}{\delta
q^j})=\frac{An}{2}G_{ij},
\\ \mbox{ } \\
 Ric(\frac{\partial}{\partial
p_i},\frac{\partial}{\partial p_j})=\frac{An}{2}H^{ij},
\\ \mbox{ } \\
Ric(\frac{\partial}{\partial p_i},\frac{\delta}{\delta
q^j})=Ric(\frac{\delta}{\delta q^j},\frac{\partial}{\partial
p_i})=0.
\end{array}
\right.
$$
Thus
\begin{equation}
 Ric=\frac{An}{2}G.
\end{equation}

 By straightforward computation, using the
relations (16),(18) and the package Ricci, the following formulas
are obtained:
\begin{equation}
\left\{
\begin{array}{l}
\frac{\delta}{\delta
q^l}QQQ^h_{ijk}=-\Gamma^h_{ls}QQQ^s_{ijk}+\Gamma^s_{li}QQQ^h_{sjk}+\Gamma^s_{lj}QQQ^h_{isk}+\Gamma^s_{lk}QQQ^h_{ijs},
\\ \mbox{ } \\
\frac{\delta}{\delta q^l}PPQ^{ijh}_k=\ \
\Gamma^s_{lk}PPQ^{ijh}_s-\Gamma^i_{ls}PPQ^{sjh}_k-\Gamma^j_{ls}PPQ^{ish}_k-\Gamma^h_{ls}PPQ^{ijs}_k,
\\ \mbox{ } \\
\frac{\delta}{\delta
q^l}PQQ^i_{jkh}=-\Gamma^i_{ls}PQQ^s_{jkh}+\Gamma^s_{lj}PQQ^i_{skh}+\Gamma^s_{lk}PQQ^i_{jsh}+\Gamma^s_{lh}PQQ^i_{jks},
\\ \mbox{ } \\
\frac{\delta}{\delta q^l}PQP^{ikh}_j=\ \
\Gamma^s_{lj}PQP^{ikh}_s-\Gamma^i_{ls}PQP^{skh}_j-\Gamma^k_{ls}PQP^{ish}_j-\Gamma^h_{ls}PQP^{iks}_j,
\\ \mbox{ } \\
\frac{\partial}{\partial
p_l}QQQ^h_{ijk}=-P^{hl}_sQQQ^s_{ijk}+P^{sl}_iQQQ^h_{sjk}+P^{sl}_jQQQ^h_{isk}+P^{sl}_kQQQ^h_{ijs},
\\ \mbox{ } \\
\frac{\partial}{\partial p_l}PPQ^{ijh}_k=\ \
P^{sl}_kPPQ^{ijh}_s-P^{il}_sPPQ^{sjh}_k-P^{jl}_sPPQ^{ish}_k-P^{hl}_sPPQ^{ijs}_k,
\\ \mbox{ } \\
\frac{\partial}{\partial
p_l}PQQ^i_{jkh}=-P^{il}_sPQQ^s_{jkh}+P^{sl}_jPQQ^i_{skh}+P^{sl}_kPQQ^i_{jsh}+P^{sl}_hPQQ^i_{jks},
\\ \mbox{ } \\
\frac{\partial}{\partial p_l}PQP^{ikh}_j=\ \
P^{sl}_jPQP^{ikh}_s-P^{il}_sPQP^{skh}_j-P^{kl}_sPQP^{ish}_j-P^{hl}_sPQP^{iks}_j.
\end{array}
\right.
\end{equation}
Due to the relations (14),(17), we have
$$
(\nabla_\frac{\delta}{\delta q^l}K)(\frac{\delta}{\delta
q^i},\frac{\delta}{\delta q^j})\frac{\delta}{\delta
q^k}=(\frac{\delta}{\delta
q^l}QQQ^h_{ijk}+\Gamma^h_{ls}QQQ^s_{ijk}-\Gamma^s_{li}QQQ^h_{sjk}-
\Gamma^s_{lj}QQQ^h_{isk}-\Gamma^s_{lk}QQQ^h_{ijs})
\frac{\delta}{\delta q^h}+
$$
$$
~~~~~~~~~~~~~~~~~~~~~~~~~~~~~~~~~~~~~~~~~~~~~~+(S_{hls}QQQ^s_{ijk}
+S_{slk}QQQ^s_{ijh}+S_{slj}PQQ^s_{ikh}
-S_{sli}PQQ^s_{jkh})\frac{\partial}{\partial p_h}.
$$

The coefficient of $\frac{\delta}{\delta q^h}$ is zero due to the
relations (20). By straightforward computation, using the
relations (16),(18) and the package Ricci, we obtain that the
coefficient of $\frac{\partial}{\partial p_h}$ is zero. Thus
$$
(\nabla_\frac{\delta}{\delta q^l}K)(\frac{\delta}{\delta
q^i},\frac{\delta}{\delta q^j})\frac{\delta}{\delta q^k}=0.
$$

Similarly
$$
(\nabla_\frac{\partial}{\partial p_l}K)(\frac{\delta}{\delta
q^i},\frac{\delta}{\delta q^j})\frac{\delta}{\delta
q^k}=(\frac{\partial}{\partial
p_l}QQQ^h_{ijk}+P^{hl}_sQQQ^s_{ijk}-P^{sl}_iQQQ^h_{sjk}-P^{sl}_jQQQ^h_{isk}
-P^{sl}_kQQQ^h_{ijs})\frac{\delta}{\delta
q^h}.
$$

The coefficient of $\frac{\delta}{\delta q^h}$ is zero due to the
relations (20). Thus
$$
(\nabla_\frac{\partial}{\partial p_l}K)(\frac{\delta}{\delta
q^i},\frac{\delta}{\delta q^j})\frac{\delta}{\delta q^k}=0.
$$

Similarly, we have computed the covariant derivatives of curvature
tensor field K in the local adapted frame $(\frac{\delta}{\delta
q^i},\frac{\partial}{\partial p_i})$ with respect to the
connection $\nabla$ and we obtained in all the cases that the
result is zero . Therefore
$$
\nabla K = 0.
$$

 Now we may state our main
result.\\

{\bf Theorem 5.} {\it If the Riemannian manifold $(M,g)$ has
positive constant sectional curvature $c$, the conditions (12) are
fulfilled and the components of the metric $G$ are given by (13)
 then $(T^*_{0A}M,G,J)$ is a locally symmetric K\"ahler Einstein
manifold.}

\vskip 1.5cm

\begin{minipage}{2.5in}
\begin{flushleft}
D.D.Poro\c sniuc\\
Department of Mathematics \\
National College "M. Eminescu" \\
Boto\c sani, Rom\^ania.\\
e-mail: danielporosniuc@lme.ro
\end{flushleft}
\end{minipage}
\end{document}